\documentclass[12pt]{amsart}

\usepackage[margin=1in,marginparwidth=0.8in, marginparsep=0.1in]{geometry}

\pdfoutput=1

\renewcommand{\baselinestretch}{1.2}

\usepackage{soul}

\usepackage{amsfonts,amssymb,latexsym,amsmath,graphicx}

\usepackage[
bookmarks=true, bookmarksopen=true,%
bookmarksdepth=3,bookmarksopenlevel=2,%
colorlinks=true,%
linkcolor=blue,%
citecolor=blue,%
filecolor=blue,%
menucolor=blue,%
urlcolor=blue]{hyperref}

\usepackage{epigraph}

\usepackage[all, knot, poly]{xy}
\xyoption{knot}

\usepackage{times}
\usepackage{mathrsfs}

\newtheorem*{theorem*}{Theorem}

\theoremstyle{definition}

\newtheorem*{definition*}{Definition}

\theoremstyle{remark} 
\newtheorem*{remark}{Remark}

\newcommand{\B}{\mathbb{B}}
\newcommand{\C}{\mathbb{C}}
\newcommand{\D}{\mathbb{D}}

\renewcommand{\L}{\mathbb{L}}

\newcommand{\R}{\mathbb{R}}

\newcommand{\Z}{\mathbb{Z}}

\newcommand{\cC}{\mathcal{C}}

\author{Vivek Shende}

\title{Microlocal category for Weinstein manifolds via h-principle}

\begin{document}

\begin{abstract}
On a Weinstein manifold, we define a constructible co/sheaf of categories on the skeleton. 
The construction works with arbitrary coefficients, and depends only on the homotopy class of 
a section of the Lagrangian Grassmannian of  the stable symplectic normal bundle. 

The definition is as follows.  Take any, possibly with high codimension, exact embedding into a 
cosphere bundle. 
Thicken to a hypersurface, and consider the Kashiwara-Schapira stack along the
thickened skeleton.  Pull back along the inclusion
of the original skeleton.  

Gromov's h-principle for contact embeddings guarantees existence and uniqueness up to isotopy of such an embedding.  
Invariance of microlocal sheaves along such isotopy is well known.  

We expect, but do not prove here, invariance of the global sections of this co/sheaf of categories under Liouville deformation.
\end{abstract}

\maketitle

\thispagestyle{empty}

Let $(W, d\lambda)$ be a compact exact symplectic manifold with convex contact boundary.  Let
$X$ be the Liouville field, i.e. $d\lambda(X, \cdot) = \lambda$.  The skeleton
of such a manifold is, by definition, the locus $\Lambda \subset W$ of points which do not
escape under the flow of $X$.  We assume that $\Lambda$ is isotropic and Whitney stratifiable.  
This includes the Weinstein manifolds in the sense of \cite{EG, cieliebakeliashberg, Eli}. 

From such a manifold, one can construct the so-called ``wrapped Fukaya category'' \cite{abouzaidseidel}; 
its objects are Lagrangians equipped with various data, its morphism spaces are certain Hamiltonian trajectories, 
and its higher structures are defined as usual in terms of pseudo-holomorphic discs.  Kontsevich \cite{Kon} conjectured
that the resulting category localizes to a cosheaf of categories on the skeleton $\Lambda$.  According
to Nadler  \cite{Nwms}, at least when $W$ can be embedded as a hypersurface in a cotangent bundle, there is a natural candidate
cosheaf, coming from the microlocal sheaf theory of Kashiwara and Schapira \cite{KS}. 

Nadler's candidate cosheaf on $\Lambda$ is obtained as follows.  For a conical Lagrangian $\L$ in a cotangent bundle
$T^*M$, 
the microlocalization of \cite{KS} can be used to construct a sheaf of categories $\mu sh_\L$ on $\L$ by sheafifying the presheaf 

$$\mu sh^{pre}_\L(\Omega) := D_{T^*M \setminus (\Omega \setminus \L)} (M) / D_{T^*M \setminus \Omega}(M) $$ 

Here, $\Omega$ is an open subset of $T^*M$.  By
$D(M)$ we mean some appropriate triangulated dg or stable $\infty$-category of sheaves on $M$;\footnote{
Our methods are largely indifferent to the precise choice of coefficients; they 
work whenever the above makes sense, e.g. over a field or over spectra.  We shall therefore not be precise on this point, 
instead referring to \cite{Nwms}, and \cite[Appendix A]{Lur2} and \cite{RS, JT} for further discussions.
} by $D_X (M)$ we mean
to require the sheaves to be microsupported in $X$, in the sense of \cite{KS}.   

For foundational material regarding sheaves of categories, see \cite{Lur1, Lur2, GR}.  
According to \cite{Nwms}, the restriction maps of $\mu sh$ have both adjoints; passing to the left adjoints turns the sheaf into a cosheaf whose corestrictions
preserve compact objects.  Due to this dual nature, we refer to $\mu sh$ as a co/sheaf. While $\mu sh$ is a priori 
defined on $T^*M$, it is in fact pushed forward from $\L$.

\begin{remark}
A great virtue of $\mu sh_\L$, which may not be apparent from the above discussion, 
is its computability.  See, e.g., \cite{FLTZ, STZ, STWZ, STW, Nlg, Nwms, Ku, GS}. 
\end{remark} 
 
Returning to the case at hand, we will say that a map from an exact symplectic manifold
$(W, \lambda)$ into a contact manifold $(\mathcal{V}, \xi)$ is {\em exact} if some contact form for $\xi$ pulls
back to $\lambda$.  Such a map must be an immersion; we term it an {\em exact embedding} if, in addition, it is 
injective.  This notion is discussed in detail in \cite{Av, Eli}.  Note that 
varying the choice of the contact form, i.e. rescaling $\lambda$,  does not affect the flow lines of the 
Liouville vector field; in particular, the skeleton of $W$ is determined by the underlying map of an exact embedding
 \cite[p.6]{Eli}. 

Suppose now given an exact embedding of $(W, \lambda)$ into a cosphere bundle $S^*M$.  Then we can take 
the conical Lagrangian formed by the positive cone on $\Lambda$ plus the zero section: 
$$\L = \R_{\ge 0} \Lambda \cup T^*_M M \subset T^*M$$
and form the co/sheaf of categories $\mu sh_\Lambda := \mu sh_\L|_\Lambda$.  

\begin{definition*}
$Sh(W) := \mu sh_\Lambda(\Lambda)$.
\end{definition*}

The above ``definition'' of $Sh(W)$ suffers two obvious defects.  
First,  an exact embedding as a hypersurface in a cosphere bundle may not exist.  Second, it is not clear 
to what extent the invariant depends on the choice of such an embedding.  

\begin{remark}
Even if $W \hookrightarrow S^*M$ is flexible, e.g. if 
$W = T^* S^1 \hookrightarrow S^* \R^2$ is a ribbon for a loose Legendrian knot, the 
category $Sh(W)$ need not vanish.  The point is that the restriction to $\Lambda$ gives
a category which should be seeing only a neighborhood of $W$, in which it is not loose.  

On the other hand, the category sees the rotation: for a Legendrian knot with 
nonzero rotation,  the category contains only periodic
objects.  We refer to \cite{Gui} for a detailed discussion.
\end{remark}

Our purpose here is to resolve 
these difficulties by appeal to Gromov's h-principle for contact embeddings; see \cite{Gromov, EM, Datta}.   Let us recall that 
an ``h-principle'' says roughly that the space of solutions to some given problem is homotopy equivalent
to the space of solutions to some linearization of the problem.  In this case, recall that a formal contact
immersion $(U, \eta) \to (V, \xi)$ is any map $i: U \to V$, and any monomorphism $TU \to i^*TV$ inducing
a monomorphism of conformally symplectic vector bundles $\eta \to i^* \xi$.  
The parametric h-principle holds, meaning that the space of actual contact immersions is homotopy equivalent to the 
space of formal contact immersions.  More to the point, when $U$ is open and the inclusion
has positive codimension, then the same holds for contact embeddings (now we should require the original
map $i:U \to V$ to be injective) \cite[12.3.1]{EM}.  Finally, when $V = \R^{2n + 1 \gg \dim U}$,  the Stiefel manifolds which classify the formal
data become arbitrarily connected.   In other words, just as any manifold admits an embedding into some sufficiently large $\R^n$,
which becomes  unique upto increasingly unique homotopy as $n \to \infty$, so too every contact manifold admits an eventually unique
embedding into $\R^{2n+1 \gg 0}$.   

We apply this to our $W$ by first taking the canonical (exact) embedding
into the contactization $W \hookrightarrow \mathcal{W}$, and composing  with a contact embedding  $\mathcal{W} \hookrightarrow \R^{2n +1}$.

\begin{remark}
Note that after this stabilization, we lose the rotation.  This may cause some cognitive dissonance: we have seen
the rotation can affect the category.  We will find it again later. 
\end{remark}

One might be tempted to use the definition above on our h-principled $W \hookrightarrow \R^{2n+1} \hookrightarrow S^* \R^{n+1}$. 
However, the resulting category of sheaves would just be zero: the skeleton of $W$ would be 
isotropic but not Legendrian, where as the microsupport of any sheaf is co-isotropic \cite{KS}.  
This is the analoguous statement in sheaf theory of the fact that subcritical items 
are generally invisible from the point of view of Floer theory.

Instead we thicken to a hypersurface embedding $$\widetilde{W} = W \rtimes \B \hookrightarrow  \R^{2n+1}$$
Here, $\B$ is some Darboux ball, $W \rtimes \B$ is a neighborhood of $W$ in the restriction to $W$ of the symplectic normal bundle
$\nu_\phi$ to the embedding $\phi: \mathcal{W} \hookrightarrow \R^{2n+1}$, and
such a thickening exists by a standard neighborhood theorem; see e.g. 
\cite{Av}.   

We write $Gr(\nu_\phi) \to \mathcal{W}$ for the Lagrangian Grassmannian of the symplectic normal bundle.  Assume it
has a section.  Fixing such a section $\sigma$ allows us to choose a Lagrangian disk bundle 
$$\widetilde{\Lambda} := \Lambda \rtimes \D \subset W \rtimes \D \subset W  \rtimes \B$$ 

Fixing an embedding $\R^{2n+1} \hookrightarrow S^* \R^{n+1}$, we write $\widetilde{\L} \subset T^* \R^{n+1}$ for the conical Lagrangian 
given by the union of the zero section and the cone over $\widetilde{\Lambda}$.  

It remains the case that $\mu sh_{\widetilde{\L}}(\widetilde{\L} \setminus \R^n) = 0$. 
However, the fact that $\widetilde{W}$ came from a lower dimensional manifold by thickening 
is only visible to $\widetilde{\L}$ through its boundaries coming from the boundary of the disk $\D$. E.g., a small neighborhood
of a point in a Legendrian knot cannot be distinguished from a small neighborhood of a point in a Legendrian interval, which in
turn is the skeleton of the thickening of an embedding of $W = point$ into $S^* \R^2$.  

Consider therefore the inclusion 
$$\Lambda =
\Lambda \times \mathbf{0} \hookrightarrow \Lambda \rtimes \D  = \widetilde{\Lambda} \hookrightarrow \widetilde{\L}$$  
Evidently it stays away from the boundaries of $\D$.  

\begin{definition*}
$\mu sh_\Lambda := \mu sh_{\widetilde{\L}}|_\Lambda$
\end{definition*}

\begin{theorem*}
The co/sheaf of categories $\mu sh_\Lambda$ depends only on the homotopy type of the section
of the Lagrangian Grassmannian of the stable symplectic normal bundle. 
\end{theorem*} 

We pause to explain what is a stable symplectic normal bundle.  Recall that by
the Hirsch-Smale h-principle, any two embeddings of a given manifold $M$ into $\R^{n \gg 0}$ are isotopic; thus 
there is a precise sense in which the normal bundle to such an embedding stabilizes, defining
a class $\nu_M \in \pi_0 Map(M, BO)$, equal to the negative of the tangent bundle. 

Gromov's $h$-principle correspondingly guarantees the symplectic normal bundles to embeddings $\phi: \mathcal{W} \to \R^{2n+1}$
 stabilize to some element of $\nu_\mathcal{W} \in \pi_0 Map(\mathcal{W}, BU)$, equal to the negative
of the contact distribution.   In our construction above we needed a choice of Lagrangian sub-bundle of $\nu_\phi$.  Note that under stablization
$\R^{2n+1} \subset \R^{2n+1} \times \R^{2m}$, there is a canonical extension of this sub-bundle by just taking the product with some Lagrangian
subspace of $\R^{2m}$.  Thus we obtain a
Lagrangian sub-bundle of the stable $\nu_\mathcal{W}$, i.e. a section of the Lagrangian Grassmannian of $\nu_\mathcal{W}$. 
The assertion is that $\mu sh_\Lambda$ only depends on the homotopy class of this section. 

\begin{remark}
For a hypersurface embedding $W \hookrightarrow \R^{2n+1}$, the symplectic normal bundle
is trivial.  Upon stabilizing to $W \hookrightarrow \R^{2n+1} \times \R^{2m}$, the symplectic normal bundle is the trivial $\R^{2m}$, and 
our prescription above gives as choice of Lagrangian some fixed $\R^m$ inside.  Note however that an isotopy inside 
the large $\R^{2n+1} \times \R^{2m}$ between two such stabilized embeddings will not generally preserve this choice of section.  
Thus we recover the rotation. 
\end{remark}

\begin{proof}
We recall a standard fact in microlocal sheaf theory.   Suppose given a Legendrian $\Xi \subset S^*M$, which admits a neighborhood $N$
and a contactomorphism $$(N, \Xi) \cong (N_0, \Xi_0) \times (T^*(-1,1)^k, (-1,1)^k)$$ where 
$N_0$ is a contact manifold containing the Legendrian $\Xi_0$.  Then $\mu sh_\Xi$ is 
locally constant along the directions $(-1,1)^k$: it is pulled back from some co/sheaf on $\Xi_0$. This can be shown
either by noncharacteristic deformation arguments or using the theory of contact transformations \cite{KS}. 

To see independence of the choice of thickening $\widetilde{W}$, and the choice (within
its homotopy class) of the Lagrangian disk bundle, just consider a family connecting such choices.  

To check independence of $\phi: \mathcal{W} \hookrightarrow   \R^{2n+1}$, we again use Gromov's h-principle. 
 Suppose given another embedding, $\phi': \mathcal{W} \hookrightarrow S^* \R^{2n'+1}$, with the same
stable symplectic normal bundle as $\phi$.  For the moment let us distinguish 
$\mu sh_\Lambda := (\phi|_{\Lambda})^* \mu sh_{\widetilde{\L}}$ and $\mu sh_\Lambda ' := (\phi'|_{\Lambda})^* \mu sh_{\widetilde{\L}}$. 

By composing with inclusions $\R^{2n+1} \to \R^{2N+1}$, we may as well assume that $\phi, \phi'$ have the same codomain.  
This stabilization changes the microsupport by a trivial factor, hence does not affect $\mu sh_\Lambda$ or $\mu sh_\Lambda'$. .  
Taking $N \gg 0$ and invoking the h-principle \cite{Gromov, EM, Datta}, the embeddings $\phi, \phi'$ are isotopic through a family of embeddings.  
We carry along the chosen section of the Lagrangian Grassmannian, hence the thickening, along this isotopy.  
The Kashiwara-Schapira stack is thus locally constant along this family, hence $\mu sh_\Lambda \cong \mu sh_\Lambda'$.  
By appealing to the full strength
of the parametric $h$-principle, we learn that this isomorphism is as unique as could be desired. 
\end{proof}


\begin{remark}
Without demanding a section of the $Gr(\nu_\mathcal{W})$,  the above construction gives
a co/sheaf of categories $\widetilde{\mu sh}_\Lambda$ over $Gr(\nu_\mathcal{W})|_{\Lambda}$, locally constant in the Grassmannian direction.
The theorem is recovered by pulling back along a section.   
\end{remark}

\begin{remark}
In fact, the existence and classification of $\mu sh_\Lambda$ depends on less than a section of $Gr(\nu_\mathcal{W})$.   
Trivializing along some  $\Lambda' \subset \Lambda$ so that 
$Gr(\nu_\mathcal{W})|_{\Lambda'} \cong \Lambda' \times  U / O = \Lambda' \times Gr(\nu_{point})$, it is clear
from the construction that 
$\widetilde{\mu sh}_{\Lambda'} \cong \mu sh_{\Lambda'} \boxtimes  \widetilde{\mu sh}_{point}$. 

In other words, the twisting is in the  $\widetilde{\mu sh}_{point}$ bundle.  This is the universal Kashiwara-Schapira stack along a smooth Legendrian; 
its stalk is one's original choice of coefficient category, $\cC$.  The corresponding local 
system of categories is classified by some map $KS: U/O \to B Aut(\cC)$.  

Thus to extract a co/sheaf on $\Lambda$ from the co/sheaf $\widetilde{\mu sh}_\Lambda$ over $Gr(\nu_\mathcal{W})|_{\Lambda}$, it suffices to give a section of the $B Aut(\cC)$-bundle 
classified by 
$$\Lambda \xrightarrow{\nu_{\mathcal W}}  BU \to B(U/O) \xrightarrow{B(KS)} B^2 Aut(\cC)$$ 
The composition with $B(KS)$ can kill a lot: e.g., $Aut(D(\Z)) = \Z \times B (\Z/2)$. 

 \cite{Lur3} suggests in passing that precisely this data should be required to define a Fukaya category, save
in place of   $KS$ he takes a de-looping of the J-homomorphism.   \cite{JT} promise to eventually show
$KS = B(J)$; specialized to $\Z$-coefficients, this is in \cite{Gui}.  
\end{remark}

\begin{remark}
Nadler suggests in \cite{Nwms} a construction of $\mu sh_\Lambda$ by cutting $\Lambda$ into pieces which embed as Legendrians
in contact cosphere bundles, defining the local categories, and then gluing by contact transformation.  Such a construction, if carried
out, will yield the same category as constructed here; this can be seen e.g. by simultaneously embedding all the local charts in one space.
\end{remark}

\begin{remark}
From the expected comparison to the wrapped Fukaya category, one expects the global sections $Sh(W) := \mu sh_\Lambda(\Lambda)$
to be invariant under Liouville deformation.  Such a deformation acts nontrivially on the skeleton $\Lambda$, so one must work
significantly harder for such a result.  Invariance under Weinstein deformations would follow given enough progress in Nadler's 
`arborealization' programme \cite{Narb, Ndef, Star, Eli, ENS}, or perhaps directly using the methods of \cite{Ndef}. 

Note however that it is unknown whether Liouville isotopic
Weinstein manifolds are in fact isotopic through Weinstein manifolds, or even through manifolds with stratifiable isotropic skeleta. 
It will follow from \cite{GPS1, GPS2}, and sufficient arborealization, 
that the category defined here is invariant under Liouville deformation.  In a subsequent article I will give a 
different approach to this invariance, by geometrically wrapping within microlocal sheaf theory \cite{Shen}. 
\end{remark}

\begin{remark}
Tamarkin \cite{Tam} and Tsygan \cite{Tsy} have constructed microlocal categories associated to compact symplectic manifolds.
Recall that compactification of a Weinstein manifold deforms its Fukaya category \cite{Sher, Sei}.   One may hope the categories of
 \cite{Tam, Tsy} deform $Sh(W)$, and that their equivalence with Fukaya categories
may be shown by deforming \cite{GPS1, GPS2}. 
\end{remark}

\vspace{3mm}

{\bf Acknowledgements.}  
This note was inspired by some remarks of Laura Starkston regarding the inability 
of sheaf theorists to define invariants of Weinstein manifolds.  I thank Roger Casals and Yasha Eliashberg for 
explanations regarding the h-principle, and Xin Jin for  explanations regarding
twisting of Kashiwara-Schapira stacks. 

I am partially
supported by the NSF CAREER grant DMS-1654545, and a Sloan fellowship.

\renewcommand{\baselinestretch}{1}

\end{document}